\documentstyle[11pt]{article}

\title{TRANSFINITE DIGRAPHS}
\author{A.H. Zemanian}
\date{}

\begin{document}
\newcommand{\be}{\begin{equation}}
\newcommand{\ee}{\end{equation}}
\newcommand{\la}{\leftarrow}
\newcommand{\ra}{\rightarrow}
\newcommand{\hla}{\hookleftarrow}
\newcommand{\hra}{\hookrightarrow}
\newcommand{\dv}{\vdash}
\newcommand{\vd}{\vdash}
\newcommand{\as}{\asymp}
\newcommand{\lla}{\langle}
\newcommand{\rra}{\rangle}
\newcommand{\N}{I \kern -4.5pt N}
\newcommand{\sss}{^{*}\!}

\maketitle \baselineskip21pt

{\ Abstract --- Transfinite graphs have been defined and examined in
a variety of prior works, but transfinite digraphs had not as yet
been investigated.  The present work embarks upon such a task.  As
with the ordinals, transfinite digraphs appear in a hierarchy of
ranks indexed by the countable ordinals.  The digraphs of rank 0 are
the conventional digraphs.  Those of rank 1 are constructed by
defining certain extremities of 0-ranked digraphs, and then
partitioning those extremities to obtain vertices of rank 1.  Then,
digraphs of rank 0 are connected together at those vertices of rank
1 to obtain a digraph of rank 1.  This process can be continued
through the natural-number ranks.  However, to achieve a digraph
whose rank is the first infinite ordinal $\omega$ (i.e., the first
limit ordinal), a special kind of transfinite digraph, which we call
a digraph with an "arrow rank" must first be constructed in a way
different from those of natural-number rank.  Then, digraphs of
still higher ranks can be constructed in a way similar to that for
the natural-number ranked digraphs. However, just before each
limit-ordinal rank, a digraph of arrow rank must be set up.

Key Words:  Transfinite digraphs, ranks of digraphs, pristine
digraphs.}

\section{Introduction}

This is the second of a series of three works in which digraphs are
generalized in nonstandard and transfinite ways.  The first work
\cite{nd} defined and examined nonstandard digraphs.  The present
work does the same for transfinite digraphs.  The third work will
discuss digraphs that are both transfinite and nonstandard.

Transfinite digraphs appear in a hierarchy of ranks indexed by the
ordinals $\nu =0, 1, 2, \ldots ,$  $\omega, \omega +1, \ldots $
along with a special "arrow rank" $\vec{\omega}$ that appears after
all the natural-number ranks but precedes the first infinite-ordinal
rank $\omega$.  Other "arrow ranks" appear subsequently just before
the countable limit ordinal ranks.  A transfinite digraph of rank
$\nu$ is called a $\nu$-digraph.  A 0-digraph is a conventional
digraph.  Our definitions and analysis proceed through all the
natural-number ranks of digraphs recursively.  However, before
proceeding on to $\omega$-digraphs, we have to first introduce the
$\vec{\omega}$-digraphs through a quite different construction.
Having done so, we can then proceed recursively through the
countably infinite ordinal ranks $\nu=\omega, \omega +1,\ \omega +2,
\ldots$, but then must consider the arrow-ranked\\ $(\omega
+\vec{\omega})$-digraph before proceeding through the $\nu$-digraphs
where $\nu=\omega \cdot 2, \omega\cdot 2 +1, \omega\cdot 2 +2,
\ldots$, and so forth.  Since discussions of the $\nu$-digraphs
where $\nu>\omega$ are so similar to those for $\nu\leq\omega$, we
merely mention them in the last section.

Our notations and symbols herein are the same as those used in
\cite{nd}.

\section{Pristine transfinite digraphs}

In general, the higher-rank vertices of transfinite digraphs can
"embrace" vertices of lower ranks;  that is, a vertex of high rank
can contain a vertex of lower rank, which in turn can contain a
vertex of still lower rank, and so forth through a sequence of
vertices of decreasing ranks.  This leads to some complications in
the analysis of such digraphs.  Some simplicity can be achieved by
prohibiting such embraced vertices.  We shall do so in this
discussion of transfinite digraphs and will call the so-restricted
digraphs {\em pristine digraph}.  Every vertex of a pristine
transfinite digraph of, say, rank $\rho$ will contain only "intips"
and "outtips" of rank $\rho -1$ but no vertex of rank less than
$\rho$, as will be explicated below.

Actually, this prohibition is no restriction at all because any
digraph can be converted into a pristine digraph by "extracting" the
embraced vertices.  The process is the same as that for transfinite
graphs and is described in \cite[Section 1.4]{B7};  the only
difference arising for transfinite digraphs is that branches in the
"extraction paths" are replaced by pairs of oppositely directed,
parallel arcs.  The resulting pristine digraphs are simpler to
analyze but have an expanded structure due to the introduction of
the extraction paths.

\section{1-digraphs}

As was indicated in the Introduction, 0-digraphs are simply the
standard digraphs discussed in \cite[Section 2]{nd}.  On the other
hand, "transfinite digraphs" appear in a hierarchy of ranks of
transfiniteness, with the 0-digraphs having the initial---and not
transfinite---rank in that hierarchy.  The "1-digraphs" comprise the
first truly transfinite rank and are the subject of this section.

It turns out that basing the construction of transfinite digraphs of
ranks 2 and higher on "directed walks," rather than on "directed
paths," leads to a more general and less complicated structure.  The
reason for this is virtually the same as that encountered when
constructing undirected transfinite graphs---namely,
path-connectedness of ranks 1 and higher is not in general
transitive as a binary relation among transfinite vertices, whereas
walk-connectedness is transitive in this regard.  This matter is
discussed in \cite[Chapter 5]{B8} and will not be repeated here.  On
the other hand, the construction of digraphs of 1 can be based
either on conventional paths or on conventional walks (i.e., either
on paths or walks of rank 0);  the two methods are entirely
equivalent.  To be consistent with our subsequent constructions of
digraphs of ranks 2 or higher, we will use walk-connectedness rather
than the equivalent path-connectedness when first constructing
transfinite digraphs of rank 1.

Assume we have at hand a (standard) 0-digraph $D^{0}=\{A,V^{0}\}$.
The vertices in $V^{0}$ are now called {\em 0-vertices}.  Moreover,
the intips and outtips of an arc are assigned the rank of $-1$ and
are called $(-1)$-intips and $(-1)$-outtips, denoted by $s^{-1}$ and
$t^{-1}$, respectively.    Also, both kinds of tips are called
$(-1)$-ditips.

We start by defining a nontrivial directed walk of rank 0, also
called a (nontrivial) 0-diwalk.  A {\em 0-diwalk} $W^{0}$ is an
alternating sequence of 0-vertices $v_{m}^{0}$ and arcs $a_{m}$,
with that sequence having at least two 0-vertices: \be
W^{0}\;=\;\lla\ldots, v_{m-1}^{0},a_{m-1}, v_{m}^{0},
a_{m},v_{m+1}^{0},\ldots\rra  \label{3.1} \ee where $a_{m}=\lla
s_{m}^{-1},t_{m}^{-1}\rra$ with $s_{m}^{-1}\in v_{m}^{0}$ and
$t_{m}^{-1}\in v_{m+1}^{0}$, where the indices $m$ traverse a set of
consecutive integers, and wherein the 0-vertices and arcs may
repeat.\footnote{Compare this definition with that of a finite
dipath given by \cite[Equation (7.1)]{nd} in symbolic notation,
wherein the 0-vertices and arcs do not repeat  and the sequence is
finite.}  For the last property in particular it may happen that
$v_{m}^{0}=v_{k}^{0}$ for some $m\not= k$ and $a_{m}=a_{k}$ again
for some $m\not= k$.\footnote{It may happen that $k=m+1$ because
self-loops are permitted.}  Each arc $a_{m}$ is directed from
$v_{m}^{0}$ toward $v_{m+1}^{0}$.  For each $m$, we say that
$v_{m}^{0}$ and $a_{m}$ are {\em incident inward} and that
$v_{m+1}^{0}$ and $a_{m}$ are {\em incident outward}.  Thus, the
direction of incidence is taken from the perspective of the
arcs---not of the vertices. Moreover, we say that $W^{0}$ is
directed in the same direction as that of its arcs (i.e., from left
to right in (\ref{3.1}).

A {\em trivial 0-diwalk} is one having just one vertex and no arc.

We require that, if the sequence (\ref{3.1}) terminates on either
side, it terminates at a 0-vertex.  $W^{0}$ is call {\em two-ended}
or {\em finite} if it terminates on both sides;  $W^{0}$ is called
{\em one-ended} if it terminates on only one side;  $W^{0}$ is
called {\em endless} if it terminates on neither side.

We can define a more general kind of walk wherein the arc directions
need not conform, that is, $a_{m}$ could be either from $v_{m}^{0}$
to $v_{m+1}^{0}$ (because $s_{m}^{-1}\in v_{m}^{0}$ and
$t_{m}^{-1}\in v_{m+1}^{0}$) or from $v_{m+1}^{0}$ to $v_{m}^{0}$
(because $s_{m}^{-1}\in v_{m+1}^{0}$ and $t_{m}^{-1}\in v_{m}^{0}$).
In this more general case, we refer to (\ref{3.1}) as a {\em
0-semiwalk}---not as a 0-diwalk.  A 0-semiwalk is a walk in the
"underlying graph" of $D^{0}$.

$W^{0}$ is called {\em extended inward} (resp. {\em extended
outward}) if there exists an index $m_{0}$ such that all entries
within (\ref{3.1}) exist and are distinct for $m<m_{0}$ (resp.
$m>m_{0}$), that is, if those parts of \ref{3.1} are one-ended
0-dipaths.

A 0-{\em intip} $s^{0}$ (resp. a 0-{\em outtip} $t^{0}$ is a maximal
set of extended inward (resp. extended outward) 0-diwalks that are
pairwise identical for all $m<m_{0}$ (resp. $m>m_{0}$);  $m_{0}$
will depend upon the choice of the pair of 0-diwalks.  Any such
0-diwalk in $s^{0}$ (resp. in $t^{0}$)  will be called a {\em
representative} of $s^{0}$ (resp. of $t^{0}$).

Now, consider the set $T^{0}$ of all 0-intips and 0-outtips of the
0-digraph $D^{0}=\{A,V^{0}\}$.  We assume that $T^{0}$ is not empty.
Partition $T^{0}$ arbitrarily to get \be T^{0}\;=\;\cup_{i\in
I^{0}}T_{i}^{0}.   \label{3.2} \ee Here, $I^{0}$ is the index set of
the partition.  At this point we change notation by setting
$v_{i}^{1}=T_{i}^{0}$ and call each $v_{i}^{1}$ a {\em
1-vertex}.\footnote{We may think of the $v_{i}^{1}$ as the
connections between the extremities of $D^{0}$.}  Thus, each
$v_{i}^{1}$ is a set of 0-intips and/or 0-outtips in the chosen
partition of $T^{0}$.  Also, we set $V^{1}=\{v_{i}^{1};i\in
I^{0}\}$;  thus, $V^{1}$ is the resulting set of 1-vertices.

Then, the {\em 1-digraph} $D^{1}$ is defined as the triple: \be
D^{1}\;=\;\{A,V^{0},V^{1}\}   \label{3.3} \ee

The "underlying 1-graph" $G^{1}$ of $D^{1}$ is defined as
follows.\footnote{We use the same definitions and notations for
graphs as those explicated in \cite{B7} (or in \cite{B8} when
embraced nodes are disallowed).}  First, remove the directions of
the arcs. Thus, each arc $a=\lla s^{-1},t^{-1}\rra$ becomes a {\em
branch} $b=\lla t_{1}^{-1},t_{2}^{-1}\rra$.  In the event that there
exist two oppositely directed arcs incident to the same two
vertices, the resulting two branches are connected in parallel, but
we allow parallel branches (as well as parallel arcs).  $B$ will
denote the set of the resulting branches.  Also, each vertex of
$D^{0}$ now becomes a 0-node, and $V^{0}$ is replaced by a set
$X^{0}$ of 0-nodes.  Next, each 0-intip $s^{0}$ or 0-outtip $t^{0}$
becomes a {\em 0-tip}.  Then, the same partitioning as before of the
set $T^{0}$ of 0-ditips to get the 1-vertices is now applied to the
set of all 0-tips to get a set of {\em 1-nodes} $x_{i}^{1}$, with
$x_{i}^{1}$ being the $i$th set of the partition.  $X^{1}$ denotes
the set of these 1-nodes.  Then, the {\em underlying 1-graph}
$G^{1}$ of $D^{1}$ is \be G^{1}\;=\;\{ B,X^{0},X^{1}\}  \label{3.4}
\ee

Returning to digraphs, here is some more terminology we shall use. A
one-ended or endless 0-diwalk $W^{0}$, as given by (\ref{3.1}), is
said to {\em traverse} a 0-intip $s^{0}$ (resp. a 0-outtip $t^{0}$)
if there exists an index $m_{0}$ such that $W^{0}$ is identical to a
representative of $s^{0}$ for all $m<m_{0}$ (resp. a representative
of $t^{0}$ for all $m>m_{0}$).

Let $u^{1}$ (resp. $v^{1}$) be a 1-vertex containing $s^{0}$ (resp.
$t^{0}$), and let $W^{0}$ traverse $s^{0}$ (resp. $t^{0}$).  We then
say that $u^{1}$ and $W^{0}$ are {\em incident inward} (resp.
$v^{1}$ and $W^{0}$ are {\em incident outward}).  In this case, we
also say that $W^{0}$ {\em reaches} $u^{1}$ (resp. {\em reaches}
$v^{1}$).

At this point, we need to define a (nontrivial) 1-diwalk
$W^{1}$.\footnote{ Now, a 1-diwalk and a 1-dipath are not in general
equivalent.} This is an alternating sequence of 1-vertices and
0-diwalks: \be
W^{1}\;=\;\lla\ldots,v_{m-1}^{1},W_{m-1}^{0},v_{m}^{1},W_{m}^{0},v_{m+1}^{1}\ldots\rra
\label{3.5} \ee where the indices $m$ traverse a set of consecutive
integers and the following condition holds for each $m$:  Each
$W_{m}^{0}$ is a 0-diwalk reaching $v_{m}^{1}$ and $v_{m+1}^{1}$
(thus, $v_{m}^{1}$ and $W_{m}^{0}$ are incident inward, and
$v_{m+1}^{1}$ and $W_{m}^{0}$ are incident outward).  The elements
of (\ref{3.5}) may repeat.  The direction of $W^{1}$ is the same as
the direction of its 0-walks;  that is, from left to right in
(\ref{3.5}).  We allow a 1-diwalk to terminate on either side.  For
instance, the sequence (\ref{3.5}) may terminate on the left and/or
on the right at a 1-vertex.  Another way is for it to terminate at a
0-vertex.  For example, the left-most 0-diwalk, say, $W_{m}^{0}$ in
(\ref{3.5}) may be a one-ended 0-diwalk starting at a 0-vertex
$v_{m}^{0}$ and extending outward to reach $v_{m+1}^{1}$;  we can
represent this case by writing (\ref{3.5}) as
\[ W^{1}\;=\;\lla v_{m}^{0},W_{m}^{0},v_{m+1}^{1},W_{m+1}^{0},v_{m+2}^{1},\ldots\rra \]
Alternatively, we may have
\[ W^{1}\;=\; \lla \ldots,v_{m-2}^{1},W_{m-2}^{0},v_{m-1}^{1},W_{m-1}^{0},v_{m}^{0}\rra \]
where now $W_{m-1}^{0}$ terminates on the right at the 0-vertex
$v_{m}^{0}$ and is incident inward at $v_{m-1}^{1}$.  If $W^{1}$
terminates on both the left and the right, we call it a {\em finite}
diwalk. Finally, as a special case, we may have that the 1-diwalk is
in fact a 0-diwalk of the form \be W^{1}\;=\;\lla
v_{m}^{0},W_{m}^{0},v_{m+1}^{0}\rra  \label{3.6} \ee where now
$W_{m}^{0}$ is a finite 0-diwalk terminating at $v_{m}^{0}$ and
$v_{m+1}^{0}$.

A {\em trivial 1-diwalk} has just one element, a 1-vertex or a
0-vertex.

As an example of a nontrivial 1-digraph consider Figure 1 (located
at the end of this report). It consists of a single one-ended
1-diwalk:
\[ W^{1}\,=\,\lla v_{1}^{1},W_{1}^{0},v_{2}^{1},W_{2}^{0},v_{3}^{1},\ldots\rra \]
The first 0-diwalk $W_{1}^{0}$ in $W^{1}$ is incident inward at
$v_{1}^{1}$, proceeds downward to pass through the 0-vertex
$v_{1}^{0}$, then upward to become incident outward at the 1-vertex
$v_{2}^{1}$.  The second 0-diwalk $W_{2}^{0}$ proceeds in the same
way, being incident inward at $v_{2}^{1}$, going through
$v_{2}^{0}$, and then becoming incident outward at $v_{3}^{1}$. This
continues infinitely to traverse all of $W^{1}$.

The definition of a 1-semiwalk is obtained from the definition of a
1-diwalk (\ref{3.5}) by relaxing the restrictions on the directions
of the 0-diwalks $W_{m}^{0}$ and of the arcs in $W_{m}^{0}$.  In
other words, a 1-semiwalk in the 1-digraph $D^{1}$ corresponds to a
1-walk in the underlying 1-graph $G^{1}$ of $D^{1}$.  Again, as a
special case, a 1-semiwalk may in fact be a 0-semiwalk corresponding
to (\ref{3.6}).

Two vertices $u$ and $v$ of ranks 0 or 1 are said to be {\em
strongly 1-dwconnected} (resp. {\em unilaterally 1-dwconnected},
resp. {\em weakly 1-swconnected})\footnote{The dw" herein is an
abbreviation diwalk"., and the "sw" stands for "semiwalk".} if there
exists a 1-diwalk from $u$ to $v$ and another 1-diwalk from $v$ to
$u$ (resp. a 1-diwalk from $u$ to $v$ or from $v$ to $u$ but not
necessarily both, resp. a 1-semiwalk terminating at $u$ and $v$).

A {\em 1-strong component} is a maximal set of vertices such that
every two of those vertices are strongly 1-dwconnected;  in
addition, as a special case, a single vertex may comprise by itself
a 1-strong component connected to itself by a trivial 1-diwalk.
Fact: Every vertex lies in exactly one 1-strong component.

A {\em 1-unilateral component} is a maximal set of vertices such
that every two of them are unilaterally 1-dwconnected.  (Now, a
single vertex by itself cannot be a 1-unilateral component.)  Fact:
Every vertex lies in at least one 1-unilateral component.

A {\em 1-weak component} is a maximal set of vertices such that
every two of them are weakly 1-swconnected.  (Here again, a single
vertex by itself cannot be a 1-weak component.) Fact: Every vertex
lies in exactly one 1-weak component.

Let us note again that a 1-diwalk may in fact be a 0-diwalk.  Thus,
a 1-strong component may actually  be a 0-strong component.
Similarly, a 1-lateral component or a 1-weak component may be of 0
rank.\footnote{See \cite[Section 8]{nd} with regard to these
components of 0 rank.}

\section{2-digraphs}

Just as the definition of 1-vertices is based on partitions of set
of 0-intips and 0-outtips, the definition of  "2-vertices" uses
partitions of sets of  "1-intips" and "1-outtips".  The definitions
of the latter entities require the idea of "extended" 1-walks.

A 1-walk $W^{1}$ (see (\ref{3.5})) in the 1-digraph $D^{1}$ is
called {\em extended inward} (resp. {\em extended outward}) if there
exists an index $m_{0}$ such that all elements of (\ref{3.5}) exist
and are distinct for $m<m_{0}$ (resp. $m>m_{0}$).  A {\em 1-tip}
(resp. {\em 1-outtip}) of $D^{1}$ is a maximal set of extended
inward (resp. extended outward) 1-diwalks that are pairwise
identical for all $m<m_{0}$ (resp. $m>m_{0}$), where the index
$m_{0}$ depends upon the choice of the pair of 1-diwalks.  Let
$T^{1}$ be the set of all 1-intips and 1-outtips.  We assume that
$T^{1}$ is not empty. Partition $T^{1}$ arbitrarily to get
$T^{1}=\cup_{i\in I^{1}}T^{1}_{i}$.  Here, $I^{1}$ is the index set
of the partition. Again, we change notation by setting
$v_{i}^{2}=T_{i}^{1}$.  We call each $v_{i}^{2}$ a {\em 2-vertex}.
Also, $V^{2}$ denotes the set of all 2-vertices.  Then, the {\em
2-digraph} $D^{2}$ is defined as the quadruple: \be
D^{2}\;=\;\{A,V^{0},V^{1},V^{2}\} \label{4.1} \ee

At this point, we could continue this discussion of 2-digraphs much
as we did for 1-digraphs after they were defined.  For instance, the
"underlying 2-graph" $G^{2}$ of $D^{2}$ could be obtained by
removing the directions of the arcs to obtain a set $B$ of branches
and in turn getting sets $X_{i}$ $(i=0,1,2)$ of $i$-nodes, all of
which would lead to
\[ G^{2}\;=\;\{B,X^{0},X^{1},X^{2}\}. \]

Furthermore, we could construct in turn "2-walks", "2-intips", and
"2-outtips", a set $V^{3}$ of "3-vertices", and finally the
"3-digraph":
\[ D^{3}\;=\;\{A,V^{0},V^{1},V^{2},V^{3}\}.  \]

Moreover, still more generally we could proceed recursively to
define "$\mu$-digraphs " for all natural numbers $\mu\in\N$.  This
we do in the next section.

\section{$\mu$-digraphs}

We turn now to the construction of a transfinite $\mu$-digraph for
the natural number $\mu$ by means of recursion.  So far, we have
done so for $\mu=0,1,2$.  We now assume, that for each natural
number $\rho=0,1,\ldots,\mu-1$, the $\rho$-digraphs \be
D^{\rho}\;=\;\{A,V^{0},\ldots,V^{\rho}\}  \label{5.1} \ee have been
defined along with their associated structures such as
$\rho$-vertices and $\rho$-diwalks $W^{\rho}$.  In particular,  we
have \be W^{\rho}\;=\;\lla
\ldots,v^{\rho}_{m-1},W^{\rho-1}_{m-1},v^{\rho}_{m},W^{\rho-1}_{m},v^{\rho}_{m+1},\ldots\rra
\label{5.2} \ee where the indices $m$ traverse a consecutive set of
integers.  Here, for $\rho=0$, $W_{m}^{-1}$ is replaced by an arc
$a_{m}=\lla s_{m}^{-1},t_{m}^{-1}\rra$, and $v_{m}^{0}$ is a
0-vertex, as in (\ref{3.1}).  For $\rho=1,\ldots,\mu-1$, all the
$(\rho-1)$-diwalks $W_{m}^{\rho-1}$ are incident inward at
$v_{m}^{\rho}$ and incident outward at $v_{m+1}^{\rho}$, and thus
the direction of $W^{\rho}$ is from left to right in (\ref{5.2}).
Because $W^{\rho}$ is a diwalk, the elements in (\ref{5.2}) may
repeat.  The sequence $W^{\rho}$ may terminate on either side at an
$\alpha$-vertex $(0\leq \alpha\leq \rho)$.

All the ideas in the preceding paragraph and some others that we
have yet to discuss in this section have been defined for $\rho=0$
and $\rho=1$ in Sections 3 and 4.  For $\rho=2,\ldots,\mu$, they
will be precisely specified when we complete this cycle of our
recursive development from $\rho=\mu-1$ to $\rho=\mu$.

 The $(\mu-1)$-diwalk $W^{\mu-1}$ is call {\em extended inward} (resp. {\em extended outward})
 if there exists an index $m_{0}$ such that all
 elements within (\ref{5.2}) exist and are distinct for $m<m_{0}$ (resp. $m>m_{0}$.)
 A $(\mu-1)$-{\em intip} $s^{\mu-1}$ (resp. a $(\mu-1)$-outtip $t^{\mu-1}$) is a maximal set of
 extended inward (resp. extended outward) $(\mu-1)$-diwalks that are pairwise
 identical for all $m< m_{0}$ (resp. $m>m_{0}$);  $m_{0}$ will depend upon
 the choice of the pair of $(\mu-1)$-diwalks.  Any such $(\mu-1)$-diwalk in $s^{\mu-1}$
 (resp. $t^{\mu-1}$) will be called a {\em representative} of $s^{\mu-1}$ (resp. $t^{\mu-1}$).
 We may on occasion refer to a 1-intip or a 1-outip as a{\em 1-ditip}.

Let $T^{\mu-1}$ be the set of all $(\mu-1)$-intips and
$(\mu-1)$-outtips of the $(\mu-1)$-digraph $D^{\mu-1}$.  Here too,
we assume that $T^{\mu-1}$ is not empty.   We partition $T^{\mu-1}$
arbitrarily to get $T^{\mu-1}=\cup_{i\in I^{\mu-1}}T_{i}^{\mu-1}$,
where $I^{\mu-1}$ is the index set of the partition.  Once again, we
change notation by setting $v_{i}^{\mu}=T_{i}^{\mu-1}$ and call each
$v_{i}^{\mu}$ a {\em $\mu$-vertex}.  Thus, each $v_{i}^{\mu}$ is a
set of $(\mu-1)$-intips and/or $(\mu-1)$-outtips.  The set
$V^{\mu}=\{v_{i}^{\mu}:i\in I^{\mu-1}\}$ denotes the set of
$\mu$-vertices.

All this yields the definition of a $\mu$-digraph, given by
(\ref{5.1}) when $\rho$ is replaced by $\mu$.  In particular, for
$\rho=\mu-1$, this cycle of our recursion has yielded one more
higher rank $\mu$ of the transfinite digraphs, namely, \be
D^{\mu}\;=\;\{A,V^{0},\ldots, V^{\mu}\}.  \label{5.3} \ee

Upon repeating our prior definitions and constructions but this time
for the rank $\mu$, we obtain the  (nontrivial) $\mu$-diwalk: \be
W^{\mu}\;=\;\lla\ldots,v_{m-1}^{\mu},W_{m-1}^{\mu-1},v_{m}^{\mu},W_{m}^{\mu-1},v_{m+1}^{\mu},\ldots\rra
\label{5.4} \ee Here, each $W_{m}^{\mu-1}$ is an endless
$(\mu-1)$-diwalk that is {\em incident inward} at $v_{m}^{\mu}$
(resp. {\em incident outward} at $v_{m+1}^{\mu}$);  in other words,
$W_{m}^{\mu-1}$ reaches $v_{m}^{\mu}$ (resp. $v_{m+1}^{\mu}$) along
a representative of a $(\mu-1)$-intip $s^{\mu-1}$ in $v_{m}^{\mu}$
(resp. along a representative of a $(\mu-1)$-outtip $t^{\mu-1}$ in
$v_{m+1}^{\mu}$). As before, the direction of $W^{\mu}$ conforms
with the directions of its $W_{m}^{\mu-1}$, that is, from left to
right in (\ref{5.4}).

Here, too, the $\mu$-diwalk may terminate on the left and/or on the
right.  That is, $W^{\mu}$ may be of the form:
\[ W^{\mu}\;=\;\lla v_{m}^{\alpha},W_{m}^{\mu-1}, v_{m+1}^{\mu},W_{m+1}^{\mu-1},v_{m+2}^{\mu},\ldots\rra  \]
where $0\leq\alpha\leq\mu$ and $W_{m}^{\mu-1}$ is incident inward at
$v_{m}^{\alpha}$ along an $(\alpha-1)$-intip.  Similarly, we may
also have
\[ W^{\mu}\;=\;\lla\ldots, v_{n-1}^{\mu},W_{n-1}^{\mu-1},v_{n}^{\mu},W_{n}^{\mu-1},v_{n+1}^{\beta}\rra \]
where $0\leq\beta\leq\mu$ and $W_{n}^{\mu-1}$ is incident outward at
$v_{n+1}^{\beta}$ along a $(\beta-1)$-outtip.  Also possible as a
special case is the three-element diwalk:
\[ W^{\mu}\;=\;\lla v_{1}^{\alpha},W_{1}^{\gamma},v_{2}^{\beta}\rra \]
where $\gamma\geq\max(\alpha-1,\beta-1)$.  Thus, diwalks of ranks
lower than $\mu$ are taken to be special cases of $\mu$-diwalks.

A {\em trivial $\mu$-diwalk} has just one element, an
$\alpha$-vertex, where $0\leq\alpha \leq\mu$.

The {\em underlying $\mu$-graph} $G^{\mu}$ of $D^{\mu}$ is obtained
by removing the directions of all the arcs and thereby the
directions of all the dipaths.  Thus, each intip and outip becomes
simply  an (undirected) tip, each vertex becomes a node, and we get
a $\mu$-graph as defined in \cite[Section 2.2]{B7}.

A $\mu$-{\em semiwalk} in $D^{\mu}$ is defined as is a $\mu$-diwalk
(\ref{5.4}) except that now directions need not conform.  That is, a
tracing along a $\mu$-semiwalk may encounter arcs in opposite
directions.  This possible nonconformity in directions also holds
for all the semiwalks of ranks less than $\mu$ within $D^{\mu}$.  To
state this another way, a $\mu$-semiwalk in $D^{\mu}$ corresponds to
a $\mu$-walk in the underlying $\mu$-graph $G^{\mu}$ of $D^{\mu}$.
As was true with diwalks, a semiwalk of rank less than $\mu$ is
taken to be a special case of a $\mu$-semiwalk.

Two vertices $u$ and $v$ of ranks $\alpha$ and $\beta$
($0\leq\alpha,\beta\leq\mu$) are said to be {\em strongly
$\mu$-dwconnected} (resp. {\em unilaterally $\mu$-dwconnected},
resp. {\em weakly $\mu$-swconnected}) if there exists a $\mu$-diwalk
from $u$ to $v$ and another $\mu$-diwalk from $v$ to $u$ (resp. a
$\mu$-diwalk from $u$ to $v$ or from $v$ to $u$ but not necessarily
both, resp. a $\mu$-semiwalk terminating at $u$ and $v$).  Here,
too, the stated $\mu$-diwalk or $\mu$-semiwalk may in fact have a
rank less than $\mu$.

A {\em $\mu$-strong component} is a maximal set of vertices such
that every two of them are strongly $\mu$-dwconnected. As a special
case, a single vertex may be comprise such a maximal set by itself.
Fact: Every vertex lies in exactly on $\mu$-strong component. Note:
Since a $\mu$-diwalk may in fact be a diwalk of lower rank, a
$\mu$-strong component may actually be a $\gamma$-strong component
where $0\leq\gamma <\mu$.

A {\em $\mu$-unilateral component} is a maximal set of vertices such
that every two of them are unilaterally $\mu$-dwconnected.  The
special case of a single vertex comprising a $\mu$-unilateral
component does not arise now.  Fact: Every vertex lies in at least
one $\mu$-unilateral component.  Note: Again a $\mu$-unilateral
component may in fact be of lower rank.

A {\em $\mu$-weak component} is a maximal set of vertices such that
every two of then are weakly $\mu$-swconnected.  Fact:  Every vertex
lies in exactly one $\mu$-weak component. Note: A $\mu$-weak
component may be of lower rank as well.

\section{$\vec{\omega}$-digraphs}

We now assume that our process for constructing $\mu$-digraphs can
be continued for ever-increasing ranks $\mu$ through all the natural
numbers.  This yields a digraph $D^{\vec{\omega}}$ that has vertices
for every natural -number ranks $\mu$: \be
D^{\vec{\omega}}\;=\;\{A,V^{0},V^{1},V^{2},\ldots\}    \label{6.1}
\ee $D^{\vec{\omega}}$ will be called an $\vec{\omega}$-digraph,
where $\vec{\omega}$ is viewed as a rank  that is larger than any
natural number but precedes the first transfinite ordinal $\omega$
\cite[pages 2 and 3]{B8}.  $\vec{\omega}$ is the first of the "arrow
ranks".

We define in $D^{\vec{\omega}}$ two kinds of one-ended diwalks; they
are substantially different from the one-ended diwalks we have so
far discussed and may or may not exist in $D^{\vec\omega}$. Assuming
that they do, one kind is the {\em one-ended
$\vec{\omega}$-outdiwalk}: \be W_{o}^{\vec{\omega}}\;=\;\lla
v_{0}^{\mu},W_{0}^{\mu},v_{1}^{\mu+1},W_{1}^{\mu+1},v_{2}^{\mu+2},W_{2}^{\mu+2},\ldots\rra
\label{6.2} \ee where, for each natural number $k=0,1,2,\ldots$,
$v_{k}^{\mu+k}$ is a $(\mu+k)$-vertex, and $W_{k}^{\mu+k}$ is a
one-ended, extended outward $(\mu+k)$-diwalk that has
$v_{k}^{\mu+k}$ as its terminal vertex, is directed rightward, and
reaches $v_{k+1}^{\mu+k+1}$ through a $(\mu+k)$-outtip, that is,
$W_{k}^{\mu+k}$ is a representative of that $(\mu+k)$-outtip.  Thus,
$W^{\vec{\omega}}$ is also directed toward the right and extends
infinitely rightward.

In contrast to (\ref{6.2}), we also have a {\em one-ended
$\vec{\omega}$-indiwalk}: \be W_{i}^{\vec{\omega}}\;=\;\lla
\ldots,W_{-3}^{\mu+2},v_{-2}^{\mu+2},W_{-2}^{\mu+1},v_{-1}^{\mu+1},W_{-1}^{\mu},v_{0}^{\mu}\rra
\label{6.3} \ee In this case, for $k=\ldots,-3,-2,-1$,
$W_{k}^{\mu-k-1}$ is a one-ended, extended inward $(\mu-k-1)$-diwalk
that has $v_{k+1}^{\mu-k-1}$ as its terminal vertex, is directed
rightward, and reaches $v_{k}^{\mu-k}$ through a $(\mu-k-1)$-intip,
that is, $W_{k}^{\mu-k-1}$ is a representative of that
$(\mu-k-1)$-intip.  Again, $W^{\vec{\omega}}$ is directed rightward
but now extends infinitely toward the left.

An {\em endless $\vec{\omega}$-diwalk} $W^{\vec{\omega}}$ is
obtained by joining these two walks at $v_{0}^{\mu}$.

$W_{o}^{\vec{\omega}}$, $W_{i}^{\vec{\omega}}$, and
$W^{\vec{\omega}}$ will themselves be called {\em extended} if their
terms are all distinct except for possibly a finite number of terms.

An example of an $\vec{\omega}$-digraph that has no extended
$\vec{\omega}$-diwalk consists of an infinite set of $\rho$-diwalks
with $\rho=0,1,2,\ldots$, with only one $\rho$-diwalk for each
$\rho$, which  share only one vertex , a 0-vertex, where they all
meet.

Of course, and example of an $\vec{\omega}$-digraph having an
extended $\vec{\omega}$-diwalk is any one of $W_{o}^{\vec{\omega}}$,
$W_{i}^{\vec{\omega}}$, or $W^{\vec{\omega}}$.

An $\vec{\omega}$-outtip is a maximal set of one-ended extended
$\vec{\omega}$-outdiwalks that are pairwise eventually identical; by
this we mean that the two one-ended $\vec{\omega}$ diwalks (as given
by (\ref{6.2})) are identical except for a finite number of terms.
Similarly, an $\vec{\omega}$-intip is a maximal set on one-ended,
extended $\vec{\omega}$-indiwalks that are pairwise eventually
identical (i.e., identical except for a finite number of terms in
any two sequences like (\ref{6.3})).  Any diwalk in an
$\vec{\omega}$-outtip is called a {\em representative} of that
$\vec{\omega}$-outtip, and similarly for a diwalk in an
$\vec{\omega}$-intip.  We will refer to both the
$\vec{\omega}$-outtips and $\vec{\omega}$-intips simply as
$\vec{\omega}$-ditips.  We will use these $\vec{\omega}$-ditips to
construct $\omega$-vertices in the next section.

\section{$\omega$-digraphs}

The next higher rank for transfinite digraphs is the rank $\omega$,
where $\omega$ is the first infinite ordinal \cite[Section 1.2]{B8}.
Assuming now that we have at hand an $\vec{\omega}$-digraph
$D^{\vec{\omega}}$ having at least one $\vec{\omega}$-ditip, let
$T^{\vec{\omega}}$ be the set of all its $\vec{\omega}$-ditips.
Partition $T^{\vec{\omega}}$ arbitrarily to get
$T^{\vec{\omega}}=\cup_{i\in I^{\vec{\omega}}}
T_{i}^{\vec{\omega}}$, where $I^{\vec{\omega}}$ is the index set of
the partition.  Here, too, we change notation by setting
$v_{i}^{\omega}=T_{i}^{\vec{\omega}}$  and call $v_{i}^{\omega}$ an
$\omega$-{\em vertex}. Thus, $v_{i}^{\omega}$ is a set of
$\vec{\omega}$-ditips.  The set $V^{\omega}=\{v_{i}^{\omega}: i\in
I^{\vec{\omega}}\}$ denotes the set of all the $\omega$-vertices.

We can now state the definition of an $\omega$-{\em digraph} as
follows: \be D^{\omega}\;=\;\{A,V^{0},V^{1},\ldots, V^{\omega}\}
\label{7.1} \ee where the ellipses indicate that vertex sets for all
the natural numbers are included.

Then, we can define an $\omega$-{\em diwalk} $W^{\omega}$ in
$D^{\omega}$ as the sequence: \be W^{\omega} \;=\;\lla
\ldots,v_{m-1}^{\omega},W_{m-1}^{\vec{\omega}},v_{m}^{\omega},W_{m}^{\vec{\omega}},v_{m+1}^{\omega},\ldots\rra
\label{7.2} \ee just as we did for a $\mu$-diwalk in (\ref{5.4})
except that $v_{m}^{\omega}$ here replaces $v_{m}^{\mu}$ there and
$W_{m}^{\vec{\omega}}$ here replaces $W_{m}^{\mu-1}$ there.  In
particular, if $v_{m}^{\omega}$ and $v_{m-1}^{\omega}$ both exist in
(\ref{7.2}), then $W_{m}^{\vec{\omega}}$ is an endless
$\vec{\omega}$-diwalk that reaches $v_{m}^{\omega}$ through an
$\vec{\omega}$-intip and reaches $v_{m+1}^{\omega}$ through an
$\vec{\omega}$-outtip.

In fact, all of the discussion of Section 5 can be transferred to
this section so long as the natural-number rank $\mu-1$ is replaces
by the arrow rank $\vec{\omega}$ and $\mu$ is replaced by $\omega$.
Because of this, we will not pursue this discussion of
$\omega$-digraphs any further.

\section{Transfinite digraphs of still higher order}

With an $\omega$-digraph $D^{\omega}$ in hand, we can define
$\omega$-intips and $\omega$-outtips and, assuming these exist, can
then define $(\omega+1)$-vertices, and finally
$(\omega+1)$-digraphs:
\[ D^{\omega+1}\;=\;\{A,V^{0},V^{1},\ldots,V^{\omega},V^{\omega+1}\} \]
Proceeding recursively, we can construct, for $k=1,2,\ldots$,
$(\omega+k)$-intips and $(\omega+k)$-outtips, and then
$(\omega+k+1)$-vertices and $(\omega+k+1)$-digraphs:
\[ D^{\omega+k+1}\;=\;\{A,V^{0},V^{1},\ldots, V^{\omega},V^{\omega+1},\ldots,V^{\omega+k+1}\} \]
After all the $(\omega+m)$-ranks $(m=0,1,2,\ldots)$ have been so
traversed, we come to the next arrow rank $\omega+\vec{\omega}$,
whose digraphs $D^{\omega+\vec{\omega}}$  are constructed in the
same way as that explicated in Section 6 with the proviso that
$\vec{\omega}$ is replaced by $\omega+\vec{\omega}$.  Next in line
are the $\omega\cdot 2$-digraphs whose construction mimics that of
Section 7.  We can then proceed onward through the ranks
$\omega\cdot 2+1$, $\omega\cdot 2+2,\ldots$, and indeed through
still larger ordinal ranks.

Can we proceed through all the countable ordinals?  The answer would
be "yes" if we could develop a completely general recursive
construction for an arbitrary countable ordinal rank and also for an
arbitrary arrow rank preceding a countable limit-ordinal rank.  But,
this has yet to be done.


\begin{thebibliography}{99}

\bibitem{B7} A.H. Zemanian, {\em Pristine Transfinite Graphs and Permissive Electrical Networks}, Birkhauser, Boston, 2001.

\bibitem{B8} A.H. Zemanian, {\em Graphs and Networks: Transfinite and Nonstandard}, Birkhauser, Boston, 2004.
\bibitem{nd} A.H. Zemanian, {\em Nonstandard Digraphs}, ECE Technical Report 1, Department of Electrical Engineering,
University at Stony Brook, Stony Brook, New York, April 15, 2009.

\end{thebibliography}
\end{document}